# On ASGS framework: general requirements and an example of implementation

KAMIL KULESZA, ZBIGNIEW KOTULSKI
Institute of Fundamental Technological Research, Polish Academy of Sciences
ul.Świętokrzyska 21, 00-049, Warsaw Poland, e-mail: {kkulesza, zkotulsk}@ippt.gov.pl

**Abstract:** In the paper we propose general framework for Automatic Secret Generation and Sharing (ASGS) that should be independent of underlying secret sharing scheme. ASGS allows to prevent the dealer from knowing the secret or even to eliminate him at all. Two situations are discussed. First concerns simultaneous generation and sharing of the random, prior nonexistent secret. Such a secret remains unknown until it is reconstructed. Next, we propose the framework for automatic sharing of a known secret. In this case the dealer does not know the secret and the secret owner does not know the shares. We present opportunities for joining ASGS with other extended capabilities, with special emphasize on PVSS and pre-positioned secret sharing. Finally, we illustrate framework with practical implementation.

## 1. INTRODUCTION

Everybody knows situations, where permission to trigger certain action requires approval of several selected entities. Equally important is that any other set of entities cannot trigger the action.
Secret sharing allows to split a secret into different pieces, called shares, which are given to the participants, such that only certain group (authorized set of participants) can recover the secret.
To make this requirement realistic, one should avoid situations were some of the protocol parties have dominant position. This reasoning resulted in creation of famework for ASGS.
 Secret sharing schemes (SSS) were independently invented by George Blakley [2] and Adi Shamir [13]. Many schemes have been presented since, for instance, Asmuth and Bloom [1], Brickell [5], Karin-Greene-Hellman (KGH) [8]. SSS can work in two modes:
1. Split control. In this case the secret itself is important, hence protected by distributing its pieces to different parties. This mode is, for instance, applied to protect proprietary secrets (e.g., Coca-Cola secret formula) or cryptographic keys.
2. Authentication. The content of the secret is secondary to the fact that only participants from the authorized set are able to recover it. This property allows to authenticate/identify parties taking part in the protocol. If they are able to recover the right secret, they are the chosen/right ones.
 Once secret sharing was introduced, people started to develop extended capabilities. Some of examples are: detection of cheaters and secret consistency verification (e.g. [11], [12], [14]), multi-secret threshold schemes (e.g., [11]), pre-positioned secret sharing schemes (e.g., [11]). The other class of extended capabilities focuses on anonymity, randomness and automatization for secret sharing procedures. Anonymous and random secret sharing was studied by Blundo, Giorgio Gaggia, Stinson in [3], [4]. Some of ideas in automatic secret sharing and generation originate from the same field.
Although verification capacity can protect against cheating, it usually comes at the price. This fact is related to the paradox stated by David Chaum, that no system can simultaneously provide privacy and integrity. Alternative approach proposed recently [9] seems to be promising shortcut. Nevertheless, the simplest way to stop cheating is to eliminate misbehaving parties from the protocol.
Dealer of the secret is the entity that assigns secret shares to the participants. Usually, the dealer has to know the secret in order to share it. This gives dealer advantage over ordinary secret participants. There are situations, where such advantage can lead to abuse. For instance: Often happens, that secret dealer is not the secret owner (e.g., owner hired the dealer to share the secret due to the task complexity). In this situation owner has to disclose secret to the dealer. Such knowledge allows Dealer to make use of the secret without cooperation of the set of authorized participants.
Having in mind two modes of operation for secret sharing schemes, we propose two solutions:

1. Automatic secret generation and sharing of random, prior nonexistent secret. It allows computing and sharing the secret "on the spot", when it is not predefined. This is typical situation for authentication mode. The secret is generated at random also the secret owner is can be eliminated.

2. Automatic sharing of a known secret. Using such an automatic procedure owner can easily share the secret, without any Dealer of the secret.

It addresses problem of secret owner not trusting the dealer. It can have added feature, that even secret owner knows neither secret shares, nor their distribution. The later decreases chances of owner interfering with the shared secret.

The goal of the paper is to propose general framework for two cases characterized above. The ASGS framework/paradigm is independent of underlying secret sharing scheme. Implementation details will be strongly dependent on characteristic of particular scheme. Hence in the framework part of the paper we provide only functional descriptions of procedures and algorithms. In the appendixes we give example of one such implementation.

The article has the following outline: section 2 is devoted for preliminaries, sections 3 contains description of automatic secret generation and sharing of random secret, while section 4 deals with automatic sharing of known secret. In section 5 we briefly discuss security requirements for ASGS. Next, we state Basic Property Conjecture, which describes sufficient conditions allowing to implement ASGS for any SSS. Last section of the main part of the paper contains open problems and concluding remarks. Example of practical implementation for KGH method (see [8]) is given in appendixes. Appendix A provides all preliminaries needed further in the example, while appendixes B, C correspond to sections 3,4 of the paper, respectively. Appendixes B, C are independent and can be studies separately. In appendix D we present verification procedure that checks consistency of the secret shares resulting from ASGS.

## 2. PRELIMINARIES

First, let's start from some philosophical background. Longman's "Dictionary of Contemporary English" describes secret, as "something kept hidden or known only to a few people". Still, there are few basic questions about nature of the secret, which need to be answered:

When does the secret existence begin?

Can secret exist before it is created?

Can secret existence be described by binary variable or is it fuzzy?

Can secret exist unknown to anyone; do we need at least one secret holder?

If secret is shared, how one can verify its validity upon combining the shares?

What does it mean that secret is shared or distributed?

Search for answer to the last question resulted in the development of secret sharing schemes, as described in introduction. As for preceding ones, answers were taken for granted without posing the questions. At that time such approach was natural, because the problem stated was to providing split control over known secret. This problem is well studied and quantitative answers are formulated in the language of information theory.

We consider situations different from initial one, hence we have to answer the questions for each case considered. The answers setup general framework for ASGS, independently of underlying secret sharing scheme. We state only qualitative description and hope that some day it will reach the quantitative level.

In the paper, we make use of automatic devices and procedures. We favor the approach, that in order to ease analysis and enhance security, they should be as simple as possible. It seems that possibility of successful implementation of ASGS for the given SSS is related to existence of some basic/fundamental property of that scheme. In the section 5 we state conjecture about this fact, while in appendix A, we state such property for KGH method.

In general framework at least two such devices will be needed. First is random number generator; with output strings having good statistical properties (e.g., see [10]). Second comes the accumulator, which is a dumb, automatic device that memory cannot be accessed otherwise than by predefined functions. Let's introduce two more concepts that are needed:

**Secure communication channel**. In this paper we assume that all the communication between protocol parties is done in the way that only communicating parties know the plaintext. Whenever we use command like "send", we presume that no third party can know the message contents. There is extensive literature on this subject; interested reader can consult for instance [11].

**Encapsulation.** Entities and devices taking part in the protocol can exchange information with others only via interface. Inner state of the entity (e.g. contents of memory registers) is hidden (encapsulated) and remains unknown for external observers. Encapsulation, originating in object-oriented paradigm (e.g., see [6]), is widely used in various fields of computer science.

Finally, last but not least, the notation:

Let $s_i^{(1)}$ and $s_i^{(2)}$ be secret shares in some SSS and $S$ denote the secret shared. $S \in K$, where $K$ is the secret space.

$C(U)$ denotes combiner algorithm for the given SSS that operates on the authorized set of shares $U$. $U^{(1)} = \{s_1^{(1)}, s_2^{(1)}, ..., s_d^{(1)}\}$ and $U^{(2)} = \{s_1^{(2)}, s_2^{(2)}, ..., s_n^{(2)}\}$ are two authorized set of secret shares such that $|U^{(1)}| = d$, $|U^{(2)}| = n$ and $C(U^{(1)}) = S = C(U^{(2)})$.

$U^{(1)}$ is called authorized set of primary secret shares that is used for verification of $U^{(2)}$. Set $U^{(2)}$ is called authorized set of user secret shares or, for the reasons that will become clear later, authorized set of master secret shares.

$P_i^{(n)}$ denotes share participant that was assigned to the secret share $s_i^{(n)}$ from $U^{(n)}$.

### 3. AUTOMATIC SECRET GENERATION AND SHARING

In this section we discuss automatic secret generation and sharing of random, prior nonexistent secret. As promised first we provide answers to the questions from the section 2.

In our approach, the secret existence begins, when it is generated. However, for the secret that is generated in the form of distributed shares, moment of creation comes when shares are combined for the first time. Before that moment, secret exists only in some potential (virtual) state. Nobody knows the secret, though secret shares exist, because they have never been combined. In order to assemble it, cooperation of authorized set of participants is required. In such a situation, there are only two ways to recover secret: by guess or by cooperation of participants from the authorized set. The first situation can be feasibly controlled by the size of the secret space, while the other one is the legitimate secret recovery procedure.

Once shares are combined, the secret is recovered. Recovered secret has to be checked against original secret in order to validate it. Hence, there must exist primary (template) copy of the secret. This can be seen from different perspective: authentication mode of operation for SSS should allow to identify and validate authorized set of participants, so, the template copy is required for comparison. For instance, consider opening bank vault. One copy of the secret is shared between bank employees that can open vault (the authorized set of secret participants). Second copy is programmed into the opening mechanism. When the employees input their combined shares, it can check whether they recover proper secret.

ASGS allows computing and sharing prior nonexistent secret "on the spot". This is typical situation for authentication mode. ASGS allows to prevent the dealer from knowing the secret or even to eliminate him at all. Using proposed procedure, it is also possible to design secret that remains unknown till the time it is recovered. Such secret cannot be compromised in the traditional meaning, because it does not exist until it is recovered. The secret is generated at random. This feature is important even without eliminating the owner. It makes the secret choice "owner independent"; hence decrease chances for the owner related attack. For instance, users in computer systems have strong inclination to use as the passwords character strings that have some meaning for them. The most popular choices are spouse/kids names and cars' registration numbers.

ASGS should allow automatic secret generation, such that:
a. The generated secret is random.
b. Two copies of the secret are created. Both secret copies are created in a distributed form.

c. Nobody knows the secret till the shares from the authorized set are combined.
d. Distributed secret shares can be replicated without compromising the secret.
e. Replication of the source set into the target set having different number of elements has to be supported
f. The secret shares resulting from replication have different values then the source shares.
g. ASGS supports same type of access structure as underlying SSS.

To meet above specification we introduce the following algorithms:

**Algorithm 1:** *SetGenerateM( d , n ).*

Description: *SetGenerateM* is used to generate two distributed copies of the random secret. It produces $U^{(1)}$ and $U^{(2)}$, such that $|U^{(1)}| = d$, $|U^{(2)}| = n$. It is automatically executed by the Accumulator. ∎

So far, generation of secret sets $U^{(1)}$ and $U^{(2)}$, was described. In order to make use of the secret shares they should be distributed to secret shares participants. Shares distribution is carried out via secure communication channel.

When $|U^{(1)}| = 1$, one is dealing with degenerate case, where $s_1^{(1)} = S$. It is noteworthy that, when $|U^{(1)}| > 1$, shares assignment to different participants $P_i^{(1)}$ allows to introduce extended capabilities in the secret sharing scheme. One of instances could be split control over secret verification procedure. Algorithm *SetGenerateM* allows only two authorized sets of secret shares to be created. Usually, only $U^{(2)}$ will be available for secret participants, while $U^{(1)}$ is reserved for shares verification. Often, it is required that there are more than one authorized sets of participants. On the other hand property, used in Algorithm 1 does not allow creating more than two authorized sets. The problem is: how to further share the secret further without recovering it's value?

This question can be answered by distributed replication of $U^{(2)}$ into $U^{(3)}$. Although all participants $P_i^{(2)}$ take part in the replication, they do not disclose information allowing secret recovery. Any of $P_i^{(2)}$ should obtain no information about any of $s_i^{(3)}$. Writing these properties formally:

1. $C(U^{(2)}) = C(U^{(3)}) = S$.

2. $P_i^{(2)}$ knows nothing about any of $s_i^{(3)}$

Such approach does not compromise $S$ and allows maintaining all previously discussed ASGS features.

**Algorithm 2:** *EqualSetReplicate( $U^{(2)}$ )*

Description: *EqualSetReplicate* is used to replicate distributed secret shares into the set with the same number of elements. It uses distributed elements of $U^{(2)}$ to create and distribute set $U^{(3)}$, such $|U^{(2)}| = |U^{(3)}| = n$. It is automatically executed by the Accumulator. ∎

**Algorithm 3:** *SetReplicateToBigger( $U^{(2)}$, d )*

*SetReplicateToBigger* is used to replicate distributed secret shares into the set with the bigger number of elements. It uses distributed elements of $U^{(2)}$ to create and distribute set $U^{(3)}$, such $n = |U_2| < |U_3| = d$. It is automatically executed by the Accumulator. ∎

**Algorithm 4:** *SetReplicateToSmaller( $U^{(2)}$, d )*

*SetReplicateToSmaller* is used to replicate distributed secret shares into the set with the smaller number of elements. It uses distributed elements of $U^{(2)}$ to create and distribute set $U^{(3)}$, such $n = |U_2| > |U_3| = d$. It is automatically executed by the Accumulator. ∎

Remarks:

1. To obtain many authorized sets of participants, multiple replication of $U^{(2)}$ takes place. In such instance $U^{(2)}$ is used as the master copy (template) for all $U^{(n)}$, $n \geq 3$. For this reason it is called authorized set of master secret shares.

2. In ASGS secret shares are derived automatically with help of simple devices like Accumulator. Nevertheless we recommend, that once distributed, shares should be tested for consistency. Namely, whether participants from various authorized set can recover the same secret $S$. Further discussion about verification is postponed to section 5.

Descriptions provided above looks like a "wish list", while algorithms illustrate only basic ideas and provide functional description. Reader that is not satisfied with this level of detail, is invited to read example of implementation in appendix B.

## 4. AUTOMATIC SECRET SHARING

In this section we discuss the case where the secret is known, hence more classical, than one presented in the previous section. Our contribution comes in two parts:
1. description of the method allowing automatic sharing of the known secret using simple automatic device;
2. using first result, we provide outline of the protocol that allows Owner and Dealer contribute independently to the process of secret sharing. The protocol may have added feature, that both know neither secret shares, nor their distribution. The later decreases chances of Owner interfering with the shared secret. Such solution addresses situation where parties of the protocol do not trust one another.

To address remaining questions from the section 2, we postulate, that:
   a. The resulting secret shares are random.
   b. Minimum two copies of the secret have to exist.
   c. At least one of the copies is in the distributed form.
   d. ASGS supports same type of access structure as underlying SSS.

To share secret $S$, secret one has to generate set $U^{(o)} = \{s_1^{(o)}, s_2^{(o)}, ..., s_n^{(o)}\}$, such $C(U^{(o)}) = S$.

Automatic secret sharing algorithm takes away from the Owner responsibility for proper construction of the secret shares. Using such a method Owner can easily share the secret.

**Algorithm 5:** *FastShare( S , n )*
Description: *FastShare* is the tool that provides fast and automatic sharing for a known secret.
Parties of the protocol: Owner, Accumulator (dumb automatic device).
It takes secret $S$ from the Owner and $n$ (number of secret participants). Algorithm is automatically executed by the Accumulator beyond Owner control. It returns $U^{(o)} = \{s_1^{(o)}, s_2^{(o)}, ..., s_n^{(o)}\}$ ∎

Resulting secret shares are not protected against modification by the Owner.

Next algorithm *SafeShares* confidentially shares secret $S$ using secret sharing mask $M$ provided by the Dealer. In the method the following conditions hold:
   a. Dealer does not know $S$;
   b. Owner does not know $M$;
   c. Owner does not know secret shares and their assignment to the secret participants

**Algorithm 6:** *SafeShares*
Functional description of the algorithm follows, while information flow is illustrated on the figure 1.
Parties of the protocol: Owner, Dealer, Accumulator secret participants
1. Dealer prepares mask $M$ needed to share the secret. It can be though as anonymous ***envelopes*** that will be used to hold secret shares. The envelopes are identical and indistinguishable. Their number is equal to the number of secret shares.
2. Owner shares the secret using *FastShare*, secret shares are placed in the envelopes. The envelopes are placed in the urn. Each of secret participants is assigned one randomly chosen envelope. As the result Owner knows neither distributed shares, nor their assignment to the participants. ∎

Secret shares, which are protected by the mask, cannot be combined in order to recover the secret. Once the shares are distributed by *SafeShares*, they have to be activated by the algorithm *ActivateShares*.

**Algorithm 7:** *ActivateShares*
Functional description of the algorithm follows, while information flow is illustrated on the figure 1.

Parties of the protocol: Dealer, secret participants

Dealer provides secret participants with the information allowing them to remove the mask (extract secret shares from the envelopes).

Once procedure is completed shares belonging to the participants from the authorized set can be used to recover the secret. ∎

It is interesting to note, that before secret shares are activated, their existence is only potential. To see it from different perspective: shares cannot be described by binary variable, they rather fuzzy (with coefficient set by the proper activation probability).

**Figure 1.** Algorithms *SafeShares* and *ActivateShares*

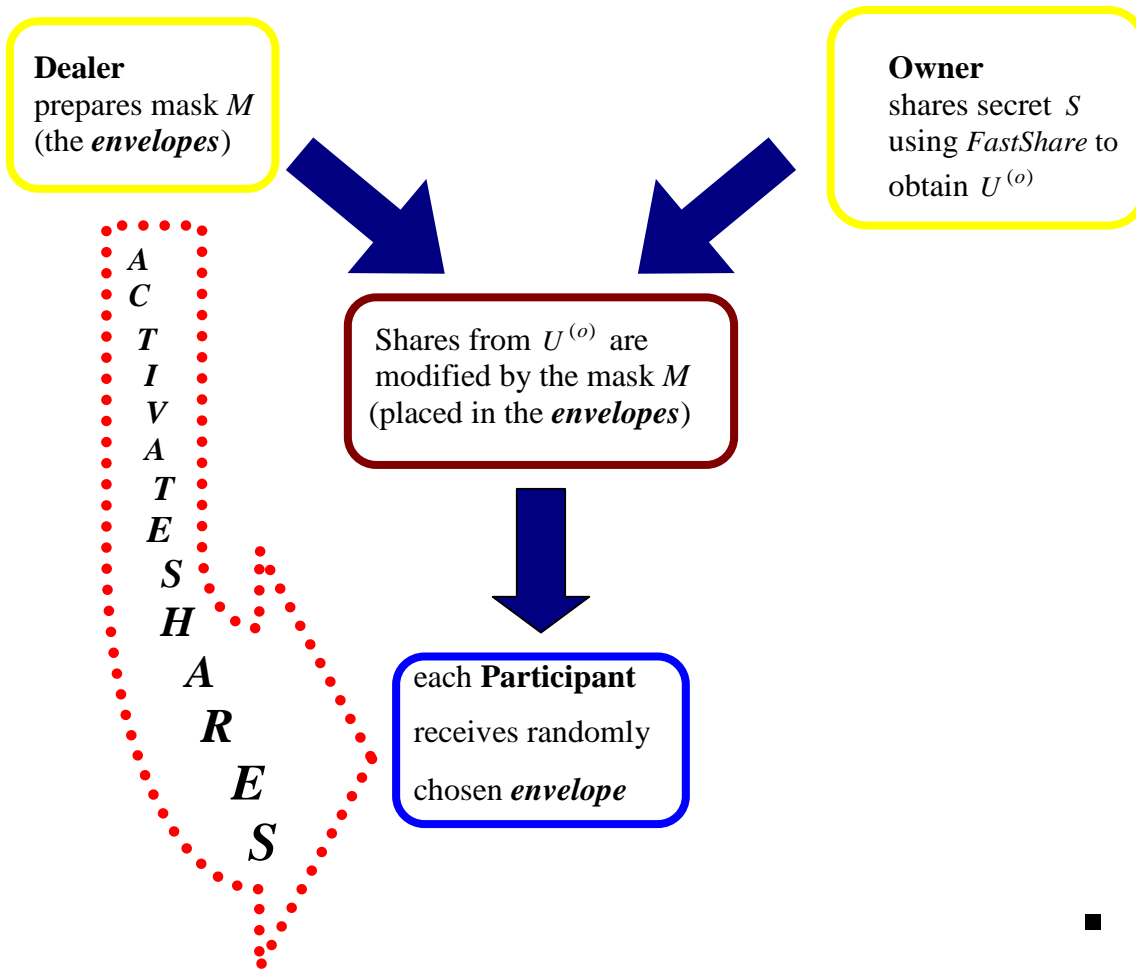

Remarks:

1. Multiple authorized sets. To create single authorized set of participants both algorithms have to be executed. Hence, to obtain many authorized sets of participants, multiple executions of *SafeShares* and *ActivateShares* take place.

2. Verification. Although much depend on the underlying SSS, proper implementation of *FastShares* should protect against cheating Owner. Yet, unless this fact is proven beyond doubt, we have to assume that Owner can modify the shares. After all, it was one of the reasons for introducing *SafeShares*. At presented level of detail one is not able to discuss cheating possibilities that might be available for the Dealer. Taking this uncertainty into account we recommend, that once activated, shares should be tested for consistency. It means, whether participants from various authorized set can recover the same secret $S$. Testing after activation allows checking both the Owner and the Dealer, as well as, whether all algorithms were properly executed. Further discussion about verification is carried out in section 5.

3. Extended capabilities. Algorithms defined above can be easily adapted to enable pre-positioned secret sharing (e.g., see [11]). In order to implement this capability it is enough to separate execution of *SafeShares* from *ActivateShares*. Hence pre-position secret sharing method can have the form:
a. Scheme is initialized by *SafeShares*.
b. It is activated using *ActivateShares*, when activation time comes,.

Again, reader that is not satisfied with this level of detail, is invited to read example of implementation in appendix C.

## 5. SECURITY RELATED STUFF

The ultimate goal is to build ASGS that will have same security features as underlying SSS. Speaking in the information theory terms: proposed framework should not decrease entropy of $S$ over $K$ (entropy of the secret over the secret space).
Qualitative description of ASGS, which was provided above, allows only making few general points concerning security of the framework. We hope, that once implementations appear, these points can be expanded into full-blown security proof. Nevertheless we make few points:
1. Security of the framework.
   Security of the framework is based on the use of secure communication channels, simple automatic devices (like Accumulator) and encapsulation principle. All these terms were specified in the section 2. Out of these three, the most unsettled are simple automatic devices. We assume, convincing security proof can be stated, as long as, devices are kept simple.
2. Basic Property Conjecture - states sufficient condition to implement ASGS for the given SSS.
   The SSS that has basic property, such that there exist operation(s) $O$ with the following characteristics:
   a. Results from performing $O$ on the authorized set of participants (possibly more then one) allow determining consistency of the secret shares.
   b. $O$ will not decrease entropy of $S$ over $K$.
   c. $O$ can be performed on the shares that are protected by *envelopes* (e.g., encrypted). ∎
   Basic property for KGH method is presented in the appendix A.
   We do not claim that only way to build ASGS for the given scheme is to find such basic property. Actually, we cannot say anything about possible alternative approach. Yet, we claim that having feasible basic property for the given SSS, we can build ASGS for that scheme.
3. Verification.
   In ASGS secret shares are derived automatically with help of the simple devices, like Accumulator. In some of the algorithms there are also interactions between parties of the protocol (e.g., Dealer provides Owner with mask $M$). Hence, we recommend, that once distributed, shares should be tested. The test has to provide information, whether participants from various authorized sets can recover the same secret $S$. The testing method depends on underlying scheme, but where possible we propose to use some Publicly Verifiable Secret Sharing (PVSS) protocol. It has to support possibility of secure testing of already distributed shares. If basic property (as described above) is found, construction of PVSS is possible. Handy example for KGH is given in appendix D.

## 6. CONCLUDING REMARKS AND FURTHER RESEARCH

The framework for ASGS was provided. Functional descriptions of the algorithms were given. Basic Property Conjecture was stated in the section 5. Finally, an example of ASGS implementation is provided in the appendixes.
A lot still needs to be done. Further research falls into three categories:
1. Research into ASGS theoretical foundations:
    - formulation of ASGS framework in terms of the information theory

- collecting more facts about Basic Property Conjecture, to state it as theorem. The final result in this field would be proof of such theorem.

2. Finding more implementations of ASGS. This includes formal proofs of security for particular implementations.

3. Placing ASGS in the broader framework within secret sharing. For instance, joining ASGS with other extended capabilities.

# APPENDIX A: Preliminaries

## KGH description

In KGH the secret is a vector of $\eta$ numbers $S_\eta = \{s_1, s_2, ..., s_\eta\}$. Any modulus $k$ is chosen, such that $k > \max(s_1, s_2, ..., s_\eta)$. All $t$ participants are given shares that are $\eta$-dimensional vectors $S_\eta^{(j)}, j = 1, 2, ..., t$ with elements in $Z_k$. To retrieve the secret they have to add the vectors component-wise in $Z_k$.

For $k = 2$, KGH method works like $\oplus$ (XOR) on $\eta$-bits numbers, much in the same way like Vernam one-time pad. If $t$ participants are needed to recover the secret, adding $t-1$ (or less) shares reveals no information about secret itself.

In practice, it is often needed that only certain specified subsets of the participants should be able to recover the secret. The authorized set of participants is a subset of all participants. Participants from such set are able to recover the secret. The access structure describes all the authorized subsets. To design the access structure with required capabilities, the cumulative array construction can be used, for details see, for example, [7], [12]. Combining cumulative arrays with KGH method, one obtains implementation of general secret sharing scheme (see, e.g., [12]).

## Remarks about procedures and algorithms presented in the appendixes.

Every routine is described in three parts:
   a. Informal description. It states the purpose of routine, describes what is being done and specifies output (when needed). Such description should be enough to comprehend the paper and get main idea behind presented methods.
   b. Routines written in pseudocode, resembling high level programming language (say C++). Level of detail is much higher than in description part. Reading through pseudocode might be tedious, but rewarding in the sense that allows appreciate proposed routines in full extend.
   c. Discussion (if needed). Methods and results are formally justified.

## Preliminaries for Algorithms

Notation:
As described in the section 2 random number generator and the Accumulator are needed. Also, secure communication channel and encapsulation have to be supported.

*RAND* yields $m_i$ obtained from a random number generator.

*ACC* denotes the value of $l$-bit memory register. Register's functions are:
*ACC.reset* sets all bits in the memory register to 0,
*ACC.read* yields *ACC*,
*ACC.store(x)* yields $ACC = ACC \oplus x$ (performs bitwise XOR of *ACC* with the input binary vector *x*, result is stored to *ACC*).

The idea of automatic secret generation and sharing, for KGH method, is based on the following property of binary vectors.

**Basic property:** Let $m_i, i = 1, 2, ..., n$, such that

$$\bigoplus_{i=1}^{n} m_i = \vec{0}, \tag{*}$$

form the set $M$. For any partition of $M$ into two disjoined subsets $C_1$, $C_2$ ($C_1 \cup C_2 = M$, $C_1 \cap C_2 = \emptyset$), it holds:

$$\bigoplus_{m_i \in C_1} m_i = \bigoplus_{m_i \in C_2} m_i. \quad \blacksquare \tag{**}$$

Now we present the procedure that generates set of binary vectors $M$.

**Procedure description:** *GenerateM* creates set of $n$ binary vectors $m_i$, satisfying condition (*). Procedure is carried out by the Accumulator. The procedure returns $M = \{m_1, m_2, \ldots, m_n\}$.

**Procedure 1:** *GenerateM( n )*

*Accumulator:*
   ACC.reset;
     for $i = 1$ to $n-1$ do
        $m_i := RAND$
        ACC.store ($m_i$)
        save $m_i$
     end //for
   $m_n = ACC.read$
   save $m_n$
   **return** $M = \{m_1, m_2, \ldots, m_n\}$
end // *GenerateM*

**Discussion**: We claim that the generated set $M$ satisfies condition (*). First, statistically independent random vectors $m_i, i = 1,2,\ldots,n-1$ are generated, while $m_n = \bigoplus_{i=1}^{n-1} m_i$, so

$$\bigoplus_{i=1}^{n} m_i = \left(\bigoplus_{i=1}^{n-1} m_i\right) \oplus m_n = \left(\bigoplus_{i=1}^{n-1} m_i\right) \oplus \left(\bigoplus_{i=1}^{n-1} m_i\right) = \vec{0} . \blacksquare$$

Further in the paper whenever we make reference to set $M$, we mean set as defined above.

# APPENDIX B: Procedures and algorithms for an automatic secret generation and sharing of random, prior nonexistent secret.

**SetGenerateM description:** It creates $U^{(1)}$ and $U^{(2)}$, such that $|U^{(1)}| = d$, $|U^{(2)}| = n$. First, *GenerateM* is used to create set $M$, such $|M| = d+n$. Next, $M$ is partitioned into $U^{(1)}$ and $U^{(2)}$. The Accumulator executes algorithm automatically.

**Algorithm 1:** *SetGenerateM( d , n )*

*Accumulator:*
   GenerateM( $d + n$ )
      for $i = 1$ to $d$ do // preparing $U^{(1)}$
         $s_i^{(1)} := m_i$
         save $s_i^{(1)}$
      end //for
      for $i = d+1$ to $d+n$ do // preparing $U^{(2)}$
         $j := i - d$
         $s_j^{(2)} := m_i$
         save $s_j^{(2)}$
      end //for
   **return** $U^{(1)} = \{s_1^{(1)}, s_2^{(1)}, \ldots, s_d^{(1)}\}$, $U^{(2)} = \{s_1^{(2)}, s_2^{(2)}, \ldots, s_n^{(2)}\}$
end// *SetGenerateM*                                                                                          ■

### Authorized set replication (same cardinality sets)

The authorized set satisfies: $|U^{(2)}| = |U^{(3)}| = n$, $U^{(2)} = \{s_1^{(2)}, s_2^{(2)}, ..., s_n^{(2)}\}$, $U^{(3)} = \{s_1^{(3)}, s_2^{(3)}, ..., s_n^{(3)}\}$. Algorithm *EqualSetReplicate* replicates set $U^{(2)}$ into the set $U^{(3)}$. It makes use of the procedure *SetReplicate*.

**SetReplicate description:** *SetReplicate* takes $U^{(2)}$ and $M$, with cardinality $|M| = 2 \cdot |U^{(2)}|$. Hence $M = \{m_1, m_2, ..., m_n\}$. First, all participants $P_i^{(2)}$ are assigned corresponding vectors $m_i$. Each of them performs bitwise XOR on their secret shares and corresponding $m_i$. Operation result is sent to the Accumulator. Accumulator adds $m_{i+n}$ to form $s_i^{(3)}$, which later is sent to $P_i^{(3)}$. As the result, simultaneous creation and distribution of $U^{(3)}$ takes place.

**Procedure 3:** *SetReplicate( M , $U^{(2)}$ )*

    *Accumulator:*
    $n := |U^{(2)}|$
        for $i = 1$ to $n$
            send $m_i$ to $P_i^{(2)}$
            $\underline{P_i^{(2)}}$: $\omega_i^{(2)} := s_i^{(2)} \oplus m_i$ // $\omega$ is *l*-bit vector (local variable)
        end//for
        for $i = 1$ to $n$
            $\underline{P_i^{(2)}}$ send $\omega_i^{(2)}$ to Accumulator
            *Accumulator:* $s_i^{(3)} := \omega_i^{(2)} \oplus m_{i+n}$
            send $s_i^{(3)}$ to $P_i^{(3)}$
        end// for
end//*SetReplicate* ∎

Algorithm *EqualSetReplicate* is the final result in this section.

**EqualSetReplicate description:** *EqualSetReplicate* takes $U^{(2)}$. It uses *SetReplicate* to create and distribute set $U^{(3)}$, such $|U^{(2)}| = |U^{(3)}| = n$.

**Algorithm 2:** *EqualSetReplicate( $U^{(2)}$ )*

    *Accumulator:*
    $n := |U^{(2)}|$
    $M := GenerateM(2n)$
    SetReplicate( $M$, $U^{(2)}$ )
end// *EqualSetReplicate*

**Discussion:** We claim that *EqualSetReplicate* fulfils requirements stated in the section 3:

1. $\bigoplus_{i=1}^{n} s_i^{(3)} = \bigoplus_{i=1}^{n} \left( s_i^{(2)} \oplus m_i \oplus m_{i+n} \right) = \left( \bigoplus_{i=1}^{n} s_i^{(2)} \right) \oplus \left( \bigoplus_{i=1}^{2n} m_i \right) = \bigoplus_{i=1}^{n} s_i^{(2)}$ as requested.

2. All $s_i^{(3)}$ result from XOR of some elements from $U^{(2)}$ with random $m_i, m_{i+n}$ hence they are random numbers. ∎

### Authorized set replication (different cardinality sets)

For $|U_2| \neq |U_3|$ there are two possibilities:

**Case 1. *SetReplicateToBigger* description**: *SetReplicateToBigger* takes $d$ and $U^{(2)}$. It generates $M$, such that $|M| = d$. Next, it uses *SetReplicate* to create and distribute first $n$ elements from $U^{(3)}$. As the result participants $P_i^{(3)}$ for $i \leq n$ have their secret shares assigned. Remaining participants $P_i^{(3)}$ are assigned $m_i$ ($i > n$) not used by *SetReplicate*. As the result $U^{(3)}$, such $n = |U_2| < |U_3| = d$ is created and distributed.

**Algorithm 3:** *SetReplicateToBigger*($U^{(2)}$, $d$)

$n := |U^{(2)}|$

$M := GenerateM(d+n)$

   *Accumulator:*

   SetReplicate($M$, $U^{(2)}$) // assigns shares for participants up to $P_n^{(3)}$, it uses first $2n$ elements of $M$

      for $i = n+1$ to $d$

         $s_i^{(3)} := m_{i+n}$

         send $s_i^{(3)}$ to $P_i^{(3)}$

      end//for

end// *BiggerSetReplicate*

**Discussion:** We claim that *SetReplicateToBigger* fulfils requirements stated in the section 3:

1. First observe, that:

$$\left(\bigoplus_{i=1}^{n}(m_i \oplus m_{i+n})\right) \oplus \left(\bigoplus_{i=n+1}^{d} m_{i+n}\right) = \left(\bigoplus_{i=1}^{2n}(m_i)\right) \oplus \left(\bigoplus_{i=2n+1}^{d+n} m_i\right) = \bigoplus_{i=1}^{d+n}(m_i)$$

so,

$$\bigoplus_{i=1}^{d} s_i^{(3)} = \left\{\bigoplus_{i=1}^{n}\left(s_i^{(2)} \oplus m_i \oplus m_{i+n}\right)\right\} \oplus \left(\bigoplus_{i=n+1}^{d} m_{i+n}\right) = \left(\bigoplus_{i=1}^{n} s_i^{(2)}\right) \oplus \left(\bigoplus_{i=1}^{d+n} m_i\right) = \bigoplus_{i=1}^{n} s_i^{(2)}$$

2. For $i > n$ all $s_i^{(3)}$ are equal to random numbers $m_i$. For $i \leq n$ all $s_i^{(3)}$ result from XOR of some elements from $U^{(2)}$ with random $m_i$, hence are random numbers. ∎

**Case 2.** *SetReplicateToSmaller* description: *SetReplicateToSmaller* takes $d$ and $U^{(2)}$. It generates $M$, such that $|M| = n + d - 1$. Next, it uses *SetReplicate* code to create $n$ secret shares $s_i^{(3)}$. First $d-1$ shares are sent to corresponding participants $P_i^{(3)}$. Remaining $s_i^{(3)}$ ($i \in \{d, d+1, ..., n\}$) are combined to form $s_d^{(3)}$ that is sent to $P_d^{(3)}$. As the result $U^{(3)}$, such that $n = |U_2| > |U_3| = d$ is created and distributed.

**Algorithm 4:** *SetReplicateToSmaller*($U^{(2)}$, $d$)

$n := |U^{(2)}|$

$M := GenerateM(n+d-1)$

   *Accumulator:*

      for $i = 1$ to $n$

         send $m_i$ to $P_i^{(2)}$

         $P_i^{(2)}: \omega_i^{(2)} := s_i^{(2)} \oplus m_i$ // $\omega$ is $l$-bit vector (local variable)

      end//for

      for $i = 1$ to $d-1$

         $P_i^{(2)}$ send $\omega_i^{(2)}$ to Accumulator

         *Accumulator:* $s_i^{(3)} := \omega_i^{(2)} \oplus m_{i+n}$

         send $s_i^{(3)}$ to $P_i^{(3)}$

      end// for

   *ACC.reset*

      for $i = d$ to $n$ // all $\omega_i^{(3)}$ for $i \leq d$ were already used

         $P_i^{(2)}$ send $\omega_i^{(2)}$ to Accumulator

*Accumulator:* ACC.store($\omega_i^{(2)}$)
      end//for
   $s_d^{(3)} = ACC.read$ // $s_d^{(3)} := \bigoplus_{i=l}^{n} \omega_i^{(2)}$

   send $s_d^{(3)}$ to $P_d^{(3)}$
end// *SetReplicateToSmaller*

**Discussion**: We claim that *SetReplicateToSmaller* fulfils requirements stated in the section 3:

1. First observe, that:

$$\left\{ \bigoplus_{i=1}^{d-1}(m_i \oplus m_{i+n}) \right\} \oplus \left( \bigoplus_{i=d}^{n} m_i \right) = \left( \bigoplus_{i=1}^{n} m_i \right) \oplus \left( \bigoplus_{i=1}^{d-1} m_{i+n} \right) = \left( \bigoplus_{i=1}^{n} m_i \right) \oplus \left( \bigoplus_{i=n}^{n+d-1} m_i \right) = \bigoplus_{i=1}^{n+d-1} m_i$$

so,

$$\bigoplus_{i=1}^{d} s_i^{(3)} = \bigoplus_{i=1}^{d-1}\left( s_i^{(2)} \oplus m_i \oplus m_{i+n} \right) \oplus \left( \bigoplus_{i=d}^{n}\left( s_i^{(2)} \oplus m_i \right) \right) = \bigoplus_{i=1}^{n} s_i^{(2)} \oplus \left( \bigoplus_{i=1}^{n+d-1} m_i \right) = \bigoplus_{i=1}^{n} s_i^{(2)}$$

2. All $s_i^{(3)}$ result from XOR of some elements from $U^{(2)}$ with random $m_i$ hence they are random numbers. ∎

# APPENDIX C: Procedures and algorithms for an automatic sharing of a known secret.

*FastShare* **description**: It takes secret $S$ and $n$ (number of secret participants). Accumulator generates random $s_i^{(o)}, i = 1, 2, ..., n-1$. Every $s_i^{(o)}$ is added to the ACC and simultaneously saved. To obtain $s_n^{(o)}$ the secret $S$ is added to ACC. Next, ACC value is read and saved as $s_n^{(o)}$. Algorithm returns $U^{(o)} = \{s_1^{(o)}, s_2^{(o)}, ..., s_n^{(o)}\}$.

**Algorithm 5:** *FastShare(S, n)*

   *Accumulator:*
   ACC.reset
      for $i = 1$ to $n-1$
         $s_i^{(o)} := RAND$
         ACC.store($s_i^{(o)}$)
         save $s_i^{(o)}$
      end// for
   *Owner:* Send secret to Accumulator
   *Accumulator:*
   ACC.store(S) //adding secret to the accumulator
   $s_n^{(o)} := ACC.read$
   save $s_n^{(o)}$
   return $U^{(o)} = \{s_1^{(o)}, s_2^{(o)}, ..., s_n^{(o)}\}$
end// *FastShare*

**Discussion**:

1. We claim that *FastShare* produces random secret shares, due to the fact that all of them originate from a random number generator. First $n-1$ shares are purely random, while the last one results from bitwise XOR of the secret and random number. More formally, $s_n^{(o)} = \bigoplus_{i=1}^{n-1} s_i^{(o)} \oplus S$. So, $s_n^{(o)}$ is random.

2. All secret shares combine to $S$. Just observe: $\bigoplus_{i=1}^{n} s_i^{(o)} = \bigoplus_{i=1}^{n-1} s_i^{(o)} \oplus s_n^{(o)} = \bigoplus_{i=1}^{n-1} s_i^{(o)} \oplus \left( \bigoplus_{i=1}^{n-1} s_i^{(o)} \oplus S \right) = S$. ∎

*SafeShares* **description**: *SafeShares* requires cooperation of two parties: Dealer and Owner. First, Dealer uses *GenerateM* to create secret sharing mask $M$, such that $\bigoplus_{m_i \in M} m_i = \vec{0}$. He also creates set $K$ of encryption keys $k_i$, such that $\bigoplus_{k_i \in M} k_i \neq \vec{0}$. $M$ elements are encrypted using corresponding keys from $K$ to form encrypted mask set $C$. $\begin{bmatrix} m_1 \\ m_2 \\ \vdots \\ m_n \end{bmatrix} \oplus \begin{bmatrix} k_1 \\ k_2 \\ \vdots \\ k_n \end{bmatrix} = \begin{bmatrix} c_1 \\ c_2 \\ \vdots \\ c_n \end{bmatrix}$ or $M \oplus K = C$

Dealer stores $K$ and sends $C$ to the Owner. Owner shares original secret $S$ using *FastShare* to obtain $S^{(o)}$. Using $C$ and $S^{(o)}$ he obtains $U^{(p)}$, which elements are randomly distributed to the participants.

$\begin{bmatrix} c_1 \\ c_2 \\ \vdots \\ c_n \end{bmatrix} \oplus \begin{bmatrix} s_1^{(o)} \\ s_2^{(o)} \\ \vdots \\ s_n^{(o)} \end{bmatrix} = \begin{bmatrix} s_1^{(p)} \\ s_2^{(p)} \\ \vdots \\ s_n^{(p)} \end{bmatrix}$ or $C \oplus U^{(o)} = U^{(p)}$

Participants receive secret shares from $S^{(p)}$ and store them.

**Algorithm 6:** *SafeShares*

Dealer:
   GenerateM(n)
   *ACC.reset*
Owner:
   FastShare(S, $n$)
     for $i = 1$ to $n$
       Dealer:
         Label $<k_i$ generation$>$:
         $k_i := RAND$
         *ACC.store*($k_i$)
         if ($i == n$ AND *ACC.read* $== \vec{0}$) {
         *ACC.store*($k_i$) //remove $k_i$ from *ACC*
         go to $<k_i$ generation$>$ // generate $k_i$ again
         } // end if
         save $k_i$
         $c_i := m_i \oplus k_i$
         send $c_i$ to Owner
       Owner: $s_i^{(p)} := c_i \oplus s_i^{(o)}$
         send $s_i^{(p)}$ to randomly chosen $P_j$ [1]
      Participant $P_j$:
         $s_j^{(p)} := s_i^{(p)}$ // share index $i$ is updated
         save $s_j^{(p)}$ // participant stores his secret share[2]
     end//for
end//*SafeShares*

---

[1] $j \in \{1, 2, ..., n\}$, one participant is allowed to obtain only one secret share. Once $s_i^{(o)}$ is send to particular $P_j$, this participant is removed from the set of participants eligible to obtain secret share.
[2] Secret share $s_i^{(p)}$ that was sent to the participant $P_j$ has now the same index $j$ as the participant.

**Discussion**:

1. Note, that $\bigoplus_{i=1}^{n} s_i^{(p)} \neq S$. So, all secret participants, upon combining their shares, will not receive $S$.

2. The rest of discussion is postponed after Algorithm 7. ■

*ActivateShares* description: *ActivateShares* requires cooperation of two parties: Dealer and Owner (of the secret). Dealer contacts participant $P_1$. Once participant's identity is established participant obtains one key from the set $K$. Participant combines $k_i$ with $s_i^{(p)}$ to obtain activated share $s_i^{(a)}$. Action is repeated for all participants.

The algorithm yields $U^{(a)} = U^{(p)} \oplus K$, where $U^{(a)} = \{s_1^{(a)} \quad s_2^{(a)} \quad ... \quad s_n^{(a)}\}$.

**Algorithm 7: *ActivateShares***

    for $i = 1$ to $n$
        Dealer:
        contacts $P_i$
        starts identification procedure
        if (*identification* == 1) sends $k_i$ to $P_i$
        Participant $P_i$:
        $s_i^{(a)} := s_i^{(p)} + k_i$
        saves $s_i^{(a)}$ // activated share is stored
    end//for
end//*ActivateShares*

**Discussion**: Once secret shares are activated, $S$ can be recovered by standard KGH procedure. We claim that $\bigoplus_{i=1}^{n} s_i^{(a)} = S$. To see it one has to combine results from two algorithms *SafeShares* and *ActivateShares*:

$$\bigoplus_{i=1}^{n} s_i^{(a)} = \bigoplus_{i=1}^{n}\left(s_i^{(p)} \oplus k_i^{(d)}\right) = \bigoplus_{i=1}^{n}\left(c_i \oplus s_i^{(o)} \oplus k_i^{(d)}\right) == \bigoplus_{i=1}^{n}\left(m_i \oplus k_i^{(p)} \oplus s_i^{(o)} \oplus k_i^{(d)}\right) = \bigoplus_{i=1}^{n}\left(m_i \oplus s_i^{(o)} \oplus k_i^{(p)} \oplus k_i^{(d)}\right)$$

One should note that particular participant $P_i$ usually obtains two different keys from $K$. Key $k_i^{(p)}$ comes from the Owner embedded in $s_i^{(p)}$, while $k_i^{(d)}$ comes from the Dealer as a part of *ActivateShares*. This does not make a difference, because $\bigoplus_{i=1}^{n}\left(k_i^{(p)} \oplus k_i^{(d)}\right) = K \oplus K = \vec{0}$ and

$$\bigoplus_{i=1}^{n} s_i^{(a)} = \bigoplus_{i=1}^{n}\left(m_i \oplus s_i^{(o)}\right) = \bigoplus_{i=1}^{n} s_i^{(o)} = S \text{ (since } \bigoplus_{i=1}^{n} m_i = \vec{0}\text{)} ■$$

# APPENDIX D: Publicly Verifiable Secret Sharing (PVSS) for KGH.

We propose efficient PVSS for KGH scheme that allows to check whether shares of the secret are correctly distributed and protects against cheating dealer. Method provides verification, without compromising the secret. Using it, secret participants in various authorized sets, can verify whether upon combining they will recover the same value of the secret $S$.

PVSS is collection of algorithms that provide method for public verification of secret shares. Scheme works for the secret sharing schemes with two or more authorized sets of participants. We present the case with two authorized sets of participants.

First, dealer has to share the secret into two sets $U^{(1)} = \{s_1^{(1)}, s_2^{(1)}, ..., s_h^{(1)}\}$, $U^{(2)} = \{s_1^{(2)}, s_2^{(2)}, ..., s_g^{(2)}\}$, such that $C(U^{(1)}) = \bigoplus_{i \in U^{(1)}} s_i^{(1)} = S = \bigoplus_{i \in U^{(2)}} s_i^{(2)} = C(U^{(2)})$.

Once the secret is shared into two authorized sets, *DistributeShares&Keys* can be used.
**Algorithm description**:

*DistributeShares&Keys* uses *RAND* to obtain random encryption keys for the secret shares. For every secret share, the key value is sent to the corresponding secret participant. The value of the secret share encrypted using derived encryption key is made public. Procedure is performed for both authorized sets of participants ($U^{(1)}, U^{(2)}$).

**Algorithm 8:** *DistributeShares&Keys*

*ACC.reset*
    for $i = 1$ to $h$
        $k_i^{(1)} := RAND$
        send $k_i^{(1)}$ to $P_i^{(1)}$
        // participant $P_i^{(1)}$ obtains his key
        $c_i^{(1)} = s_i^{(1)} \oplus k_i^{(1)}$
        publish $c_i^{(1)}$
    end// for
    for $i = 1$ to $g$
        $k_i^{(2)} := RAND$
        send $k_i^{(2)}$ to $P_i^{(2)}$
        // participant $P_i^{(2)}$ obtains his key
        $c_i^{(2)} = s_i^{(2)} \oplus k_i^{(2)}$
        publish $c_i^{(2)}$
    end// for
end// *DistributeShares&Keys*

**Discussion**:
1. When algorithm is completed, the following hold:
a. The encrypted value of secret shares $c_i^{(1)}$, $c_i^{(2)}$ are publicly available for all secret shares from $U^{(1)}, U^{(2)}$;
b. The secret key $k_i^{(1)} / k_i^{(2)}$ is known only to the participant concerned.

2. Participants from the authorized set are able to recover the secret. For instance, consider participants from $U^{(1)}$, encrypted shares are formed as shown on the left-hand side below:

$$\begin{bmatrix} s_1^{(1)} \\ s_2^{(1)} \\ \vdots \\ s_h^{(1)} \end{bmatrix} \oplus \begin{bmatrix} k_1^{(1)} \\ k_2^{(1)} \\ \vdots \\ k_h^{(1)} \end{bmatrix} = \begin{bmatrix} c_1^{(1)} \\ c_2^{(1)} \\ \vdots \\ c_h^{(1)} \end{bmatrix}$$ . When it comes to secret recovery $$\begin{bmatrix} c_1^{(1)} \\ c_2^{(1)} \\ \vdots \\ c_h^{(1)} \end{bmatrix} \oplus \begin{bmatrix} k_1^{(1)} \\ k_2^{(1)} \\ \vdots \\ k_h^{(1)} \end{bmatrix} = \begin{bmatrix} s_1^{(1)} \\ s_2^{(1)} \\ \vdots \\ s_h^{(1)} \end{bmatrix}$$ , hence the secret shares can be computed, allowing to combine the secret $\bigoplus_{i \in U^{(1)}} s_i^{(1)} = S$.

When the keys and the encrypted shares are distributed, the secret participants can verify validity of their shares (or check dealer's honesty). The first step is made with the algorithm *RecoverXORedKeys*.

**Algorithm description**:

*RecoverXORedKeys* uses Accumulator to combine all secret keys $k_i^{(1)} / k_i^{(2)}$. One round: a key from $U^{(1)}$ is sent to the Accumulator, next a key from $U^{(2)}$ is sent. Operation is repeated until all the keys are in the Accumulator.

**Algorithm 9:** *RecoverXORedKeys($U^{(1)}, U^{(2)}$)*

*ACC.reset*

If $(n \geq m)$ then *counter* := $n$
else *counter* := $m$
  for $i = 1$ to *counter*
    Participant $P_i^{(1)}$:
      if $(n \geq i)$ send $s_i^{(1)}$ to Accumulator
      else $s_i = \vec{0}$ to Accumulator
    Accumulator:
      *ACC.store*( $s_i^{(1)}$ )
    Participant $P_i^{(2)}$:
      if $(m \geq i)$ send $s_i^{(2)}$ to Accumulator
      else $s_i = \vec{0}$ to Accumulator
    Accumulator:
      *ACC.store*( $s_i^{(2)}$ )
  end// for
  *RecoverXORedKeys* := *ACC.read*
end// *RecoverXORedKeys*

**Discussion**:
1. Method of sending keys to the Accumulator enforce security upon value of combined keys from one set. At no point in time, value of the combined keys cannot be recovered using Accumulator contents. ∎

Having a way to recover the key, we are ready for algorithm *Verify*.

**Algorithm 10: Verify**$(U^{(1)}, U^{(2)})$

1. Publicly XOR all encrypted secret shares,
   store result in *XOREncryptedShares*
2. Run RecoverXORedKeys, store result in
   *XORedKeys*
3. If *(XORedKeys== XOREncryptedShares)*
   verification *POSITIVE*
   else verification *NEGATIVE*
end// *Verify*

**Discussion**:
1. We claim that if dealer is honest *XORedKeys= XOREncryptedShares*. Note that:

$$\left(\bigoplus_{i \in U^{(1)}} c_i^{(1)}\right) \oplus \left(\bigoplus_{i \in U^{(2)}} c_i^{(2)}\right) = \left(\bigoplus_{i \in U^{(1)}} \left(k_i^{(1)} \oplus s_i^{(1)}\right)\right) \oplus \left(\bigoplus_{i \in U^{(2)}} \left(k_i^{(2)} \oplus s_i^{(2)}\right)\right) = \left(\bigoplus_{i \in U^{(1)}} k_i^{(1)}\right) \oplus \left(\bigoplus_{i \in U^{(2)}} k_i^{(2)}\right) \oplus$$

$$\oplus \left(\bigoplus_{i \in U^{(1)}} s_i^{(1)}\right) \oplus \left(\bigoplus_{i \in U^{(2)}} s_i^{(2)}\right) = \left(\bigoplus_{i \in U^{(1)}} k_i^{(1)}\right) \oplus \left(\bigoplus_{i \in U^{(2)}} k_i^{(2)}\right) \oplus S \oplus S = \left(\bigoplus_{i \in U^{(1)}} k_i^{(1)}\right) \oplus \left(\bigoplus_{i \in U^{(2)}} k_i^{(2)}\right)$$

The relation presented above is symmetric in the sense that it does not differentiate between keys and encrypted shares. If dealer attempts to cheat (no matter whether in keys and/or encrypted shares), the equality will not hold. It is satisfied only when both authorized sets of participants receive the same value of the secret.

2. Because no information is revealed about $\bigoplus_{i \in U^{(1)}} k_i^{(1)}$ or $\bigoplus_{i \in U^{(2)}} k_i^{(2)}$ (only they XORed value is provided), hence verification is secure. ∎

Algorithm *Verify* is the final result of this section. PVSS, described so far, works for two authorized sets of the secret participants. However, method can be adapted to work for more authorized sets of participants. In such a case, the secret keys in *DistributeShares&Keys* have to be derived for all participants.